\pdfoutput=1 %arXiv admin added this line
\documentclass[11pt]{article}%
\usepackage{amsfonts}
\usepackage{amsmath}
\usepackage{amssymb}
\usepackage{graphicx}%
\newtheorem{theorem}{Theorem}

\newtheorem{definition}[theorem]{Definition}
\newtheorem{example}[theorem]{Example}

\newtheorem{proposition}[theorem]{Proposition}
\newtheorem{remark}[theorem]{Remark}

\begin{document}

\title{Solution of boundary value and eigenvalue problems for second order elliptic
operators in the plane using pseudoanalytic formal powers}
\author{Ra\'{u}l Castillo P\'{e}rez$^{\text{1}}$, Vladislav V. Kravchenko$^{\text{2}}$
\and and
\and Rabindranath Res\'{e}ndiz V\'{a}zquez$^{\text{1}}$\\$^{\text{1}}${\small SEPI, ESIME Zacatenco, Instituto Polit\'{e}cnico
Nacional, Av. IPN S/N, }\\{\small C.P. 07738, D.F. MEXICO}\\$^{\text{2}}${\small Department of Mathematics, CINVESTAV del IPN, Unidad
Queretaro, }\\{\small Libramiento Norponiente No. 2000, Fracc. Real de Juriquilla,
Queretaro, }\\{\small Qro. C.P. 76230 MEXICO e-mail:
vkravchenko@qro.cinvestav.mx\thanks{Research was supported by CONACYT, Mexico
via the research project 50424. The first named author wishes to thank support
from the CONACYT and the National Polytechnic Institute for the possibility of
a postdoctoral stay in the Department of Mathematics of the CINVESTAV in
Queretaro, as well as from the SIBE program of the National Polytechnic
Institute, Mexico.}}}
\maketitle

\begin{abstract}
We propose a method for solving boundary value and eigenvalue problems for the
elliptic operator $D=\operatorname*{div}p\operatorname*{grad}+q$ in the plane
using pseudoanalytic function theory and in particular pseudoanalytic formal
powers. Under certain conditions on the coefficients $p$ and $q$ with the aid
of pseudoanalytic function theory a complete system of null solutions of the
operator can be constructed following a simple algorithm consisting in
recursive integration. This system of solutions is used for solving boundary
value and spectral problems for the operator $D$ in bounded simply connected
domains. We study theoretical and numerical aspects of the method.

\end{abstract}

\section{Introduction}

The main numerical techniques for solving problems related to elliptic linear
partial differential equations with variable coefficients in one way or
another involve a discretization of a domain and solution of systems of
thousands of algebraic equations. Seldom the method of separation of variables
is applied due to its natural limitations related to the requirements of a
complete agreement between the geometry of the domain and the symmetry of the
coefficients. Moreover, the method of separation of variables implies solution
of Sturm-Liouville spectral problems which is not an easy task itself.

In the present paper we propose a different method based on some old and new
results from pseudoanalytic function theory \cite{Berskniga}, \cite{APFT}. Its
applicability is not so universal as the applicability of the finite
difference method or the finite element method. First of all, it is applicable
to problems in bounded domains in the plane and up to now only for the
operator $\operatorname{div}p\operatorname{grad}+q$. Moreover, at present we
can apply the method only when the equation
\begin{equation}
\left(  \operatorname{div}p\operatorname{grad}+q\right)
u(x,y)=0\label{maineqintro}%
\end{equation}
possesses a particular solution $u_{0}$ such that the function $f=p^{1/2}%
u_{0}$ is sufficiently smooth, nonvanishing in the domain of interest and
representable in the form $f=S(s)T(t)$ where $s$ and $t$ represent an
orthogonal coordinate system. This separable form of $f$ may cause
associations with the method of separation of variables. Nevertheless this is
a completely different technique, based on different ideas and free of the
mentioned above limitations of the method of separation of variables.

The heart of the method is the construction of a complete system of solutions
for (\ref{maineqintro}) in the domain of interest, complete in the sense
explained below (see \cite[Sect. 1.3]{Colton} for related ideas and additional
details). The system of solutions is used for approximating the solution of a
boundary value problem. Due to the linearity of equation (\ref{maineqintro})
after the construction of the system of solutions the problem reduces to
approximation of boundary conditions, for which a variety of methods can be
used. Here we apply the collocation method.

The complete system of solutions is constructed in the following way. The
knowledge of a particular solution of (\ref{maineqintro}) allows us to propose
a corresponding Vekua equation \cite{Vekua} closely related to
(\ref{maineqintro}) in the sense that the real part of any of its solutions
has the form $p^{1/2}u$ where $u$ is a solution of (\ref{maineqintro}), and
vice versa given $u$ one can easily construct a corresponding solution of the
Vekua equation \cite{Krpseudoan}, \cite{APFT}. The relation between
(\ref{maineqintro}) and the Vekua equation is similar to the relation between
the Laplace equation and the Cauchy-Riemann system. L. Bers developed
\cite{Berskniga}, \cite{BersFormalPowers} a theory of so-called pseudoanalytic
formal powers. They are generalizations of the analytic powers $(z-z_{0})^{n}$
in the sense that they are solutions of the corresponding Vekua equation and
behave locally like the analytic powers. The theory of Bers includes
generalizations of Taylor series, Runge's theorem and other basic facts from
analytic function theory. Thus, under certain quite natural conditions the
system of pseudoanalytic formal powers is complete in the space of all
pseudoanalytic functions (solutions of the Vekua equation) in the same sense
as the system of powers $(z-z_{0})^{n}$ is complete in the space of analytic
functions. To construct the pseudoanalytic formal powers the knowledge of a
corresponding generating sequence is required. Recently
\cite{KrRecentDevelopments}, \cite{APFT} an algorithm for construction of
generating sequences under additional conditions on the coefficients in the
Vekua equation was proposed. This implies that when $f=p^{1/2}u_{0}$ is
representable in a separable form the complete system of formal powers for the
Vekua equation associated with (\ref{maineqintro}) can be constructed
explicitly following Bers' recursive procedure.

We investigate the efficiency of the proposed method which we call MPFP, the
Method of Pseudoanalytic Formal Powers. We show its fast convergence and
compare its accuracy with that of the finite element method. In general, we
show that in problems addmitting the explicit construction of formal powers
and hence the application of the MPFP its use is advantageous compared to
other computational techniques based on discretization of the problem.

It is worth mentioning that the MPFP is a direct generalization of the method
of harmonic polynomials for solving boundary value problems for the Laplace
equation which has been considered in dozens of works (see, e.g.,
\cite{Cannon}, \cite{Genev1984}, \cite{Hozejovski}, \cite{Suetin}). Indeed, in
a special case when $p\equiv1$, $q\equiv0$ and $u_{0}\equiv1$ the
corresponding complete system of solutions constructed by means of the MPFP
coincides with the system of harmonic polynomials $\left\{  \operatorname{Re}%
(z-z_{0})^{n},\quad\operatorname{Im}(z-z_{0})^{n}\right\}  _{n=0}^{\infty}$.

The knowledge of a complete system of solutions for an equation corresponding
to any value of a spectral parameter allows one to use it for solving
eigenvalue problems. We consider this possibility in section \ref{SectEigen}.
The numerical results are highly promising, and it is clear that in the case
of eigenvalue problems as well as for boundary value problems further work
should be done in investigation of optimal ways of application of the MPFP.
For example, for solving eigenvalue problems by means of the MPFP we used the
simplest possible idea reducing the problem to calculation of zeros of a
certain determinant obtained by evaluating the first $N$ solutions from the
constructed complete system in $N$ points on the boundary of the domain under
consideration. Meanwhile, in principle, this natural approach works there
exist other techniques offering different ways of using the available exact
solution systems (see \cite[Sect. 1.13]{Alexidze}, where similar questions are discussed).

\section{Factorization of the operator $\operatorname{div}p\operatorname{grad}%
+q$.\label{SectFactorization}}

Let $\Omega$ be a domain in $\mathbf{R}^{2}$. Throughout the whole paper we
suppose that $\Omega$ is a simply connected domain. Denote $\partial
_{\overline{z}}=\frac{1}{2}\left(  \frac{\partial}{\partial x}+i\frac
{\partial}{\partial y}\right)  $ and $\partial_{z}=\frac{1}{2}\left(
\frac{\partial}{\partial x}-i\frac{\partial}{\partial y}\right)  $. By $C$ we
denote the operator of complex conjugation.

Note that the operator $\partial_{\overline{z}}$ applied to a real valued
function $\varphi$ can be regarded as a kind of gradient, and if we know that
$\partial_{\overline{z}}\varphi=\Phi$ in a whole complex plane or in a convex
domain, where $\Phi=\Phi_{1}+i\Phi_{2}$ is a given complex valued function
such that its real part $\Phi_{1}$ and imaginary part $\Phi_{2}$ satisfy the
equation
\begin{equation}
\partial_{y}\Phi_{1}-\partial_{x}\Phi_{2}=0, \label{casirot}%
\end{equation}
then we can reconstruct $\varphi$ up to an arbitrary real constant $c$ in the
following way%
\[
\varphi(x,y)=2\left(  \int_{x_{0}}^{x}\Phi_{1}(\eta,y)d\eta+\int_{y_{0}}%
^{y}\Phi_{2}(x_{0},\xi)d\xi\right)  +c
\]
where $(x_{0},y_{0})$ is an arbitrary fixed point in the domain of interest.
Note that this formula can be easily extended to any simply connected domain
by considering the integral along an arbitrary rectifiable curve $\Gamma$
leading from $(x_{0},y_{0})$ to $(x,y)$%
\begin{equation}
\varphi(x,y)=2\left(  \int_{\Gamma}\Phi_{1}dx+\Phi_{2}dy\right)  +c.
\label{Antigr}%
\end{equation}

By $\overline{A}$ we denote the integral operator in (\ref{Antigr}):%
\[
\overline{A}[\Phi](x,y)=2\left(  \int_{x_{0}}^{x}\Phi_{1}(\eta,y)d\eta
+\int_{y_{0}}^{y}\Phi_{2}(x_{0},\xi)d\xi\right)  +c.
\]
Thus if $\Phi$ satisfies (\ref{casirot}), there exists a family of real valued
functions $\varphi$ such that $\partial_{\overline{z}}\varphi=\Phi$, given by
the formula $\varphi=\overline{A}[\Phi]$.

The following result is in the core of the method proposed in the present work.

\begin{theorem}
\cite{KrJPhys06} \label{ThFactGenSchr}Let\ $p$ and $q$ be real valued
functions, $p\in C^{2}(\Omega)$ and $p\neq0$ in $\Omega$, $u_{0}$ be a
positive particular solution of the equation
\begin{equation}
(\operatorname{div}p\operatorname{grad}+q)u=0\text{\qquad in }\Omega\text{.}
\label{maineq}%
\end{equation}
Then for any real valued continuously twice differentiable function $\varphi$
the following equality holds%
\begin{equation}
\frac{1}{4}(\operatorname{div}p\operatorname{grad}+q)\varphi=p^{1/2}\left(
\partial_{z}+\frac{f_{\overline{z}}}{f}C\right)  \left(  \partial
_{\overline{z}}-\frac{f_{\overline{z}}}{f}C\right)  p^{1/2}\varphi,
\label{mainfact}%
\end{equation}
where
\begin{equation}
f=p^{1/2}u_{0}. \label{fandu}%
\end{equation}

\end{theorem}

\begin{remark}
Let $q\equiv0$. Then $u_{0}$ can be chosen as $u_{0}\equiv1$. Hence
(\ref{mainfact}) gives us the equality
\[
\frac{1}{4}\operatorname{div}(p\operatorname{grad}\varphi)=p^{1/2}\left(
\partial_{z}+\frac{\partial_{\overline{z}}p^{1/2}}{p^{1/2}}C\right)  \left(
\partial_{\overline{z}}-\frac{\partial_{\overline{z}}p^{1/2}}{p^{1/2}%
}C\right)  (p^{1/2}\varphi).
\]

\end{remark}

Let $f$ be a real function of $x$ and $y$. Consider the Vekua equation
\begin{equation}
W_{\overline{z}}=\frac{f_{\overline{z}}}{f}\overline{W}\text{\qquad in }%
\Omega\text{.}\label{Vekuamain}%
\end{equation}
This equation plays a crucial role in all that follows and hence we will call
it the \textbf{main Vekua equation}. We notice that the operator of this
equation is precisely the second factor in (\ref{mainfact}).

Denote $W_{1}=\operatorname*{Re}W$ and $W_{2}=\operatorname{Im}W$.

\begin{theorem}
\cite{KrJPhys06}\label{ThConjugate2} Let $W=W_{1}+iW_{2}$ be a solution of
(\ref{Vekuamain}). Assume that $f=p^{1/2}u_{0}$, where $u_{0}$ is a positive
solution of (\ref{maineq}) in $\Omega$. Then $u=p^{-1/2}W_{1}$ is a solution
of (\ref{maineq}) in $\Omega$, and $v=p^{1/2}W_{2}$ is a solution of the
equation
\begin{equation}
(\operatorname*{div}\frac{1}{p}\operatorname*{grad}+q_{1})v=0\qquad\text{in
}\Omega, \label{assocmaineq}%
\end{equation}
where
\begin{equation}
q_{1}=-\frac{1}{p}\left(  \frac{q}{p}+2\left\langle \frac{\nabla p}{p}%
,\frac{\nabla u_{0}}{u_{0}}\right\rangle +2\left(  \frac{\nabla u_{0}}{u_{0}%
}\right)  ^{2}\right)  . \label{q1}%
\end{equation}

\end{theorem}

Theorem \ref{ThConjugate2} shows us that as much as real and imaginary parts
of a complex analytic function are harmonic functions, the real and imaginary
parts of a solution of the main Vekua equation (\ref{Vekuamain}) multiplied by
$p^{-1/2}$ and $p^{1/2}$ respectively are solutions of the associated elliptic
equations (\ref{maineq}) and (\ref{assocmaineq}). The following natural
question arises then. We know that given an arbitrary real valued harmonic
function in a simply connected domain, a conjugate harmonic function can be
constructed explicitly such that the obtained couple of harmonic functions
represent the real and imaginary parts of a complex analytic function. What is
the corresponding more general fact for solutions of associated elliptic
equations (\ref{maineq}) and (\ref{assocmaineq}) (which we slightly
generalizing the definition of I. N. Vekua call metaharmonic functions). The
precise result is given in the following theorem.

\begin{theorem}
\cite{KrJPhys06} \label{CorConjugate}Let $f=p^{1/2}u_{0}$, where $u_{0}$ is a
positive solution of (\ref{maineq}) in a simply connected domain $\Omega$ and
$u$ be a solution of (\ref{maineq}). Then a solution $v$ of (\ref{assocmaineq}%
) with $q_{1}$ defined by (\ref{q1}) such that $W=p^{1/2}u+ip^{-1/2}v$ is a
solution of (\ref{Vekuamain}), is constructed according to the formula%
\begin{equation}
v=u_{0}^{-1}\overline{A}(ipu_{0}^{2}\partial_{\overline{z}}(u_{0}^{-1}u)).
\label{transfDarboux}%
\end{equation}
Let $v$ be a solution of (\ref{assocmaineq}), then the corresponding solution
$u$ of (\ref{maineq}) such that $W=p^{1/2}u+ip^{-1/2}v$ is a solution of
(\ref{Vekuamain}), is constructed according to the formula%
\begin{equation}
u=-u_{0}\overline{A}(ip^{-1}u_{0}^{-2}\partial_{\overline{z}}(u_{0}v)).
\label{transfDarbouxinv}%
\end{equation}

\end{theorem}

\begin{remark}
When $p\equiv1$, $q\equiv0$ and $u_{0}\equiv1$, equalities
(\ref{transfDarboux}) and (\ref{transfDarbouxinv}) turn into the well known
formulas in complex analysis for constructing conjugate harmonic functions.
\end{remark}

\section{Formal powers}

Briefly speaking formal powers are solutions of a Vekua equation
\begin{equation}
W_{\overline{z}}=aW+b\overline{W} \label{VekuaGeneral}%
\end{equation}
(with $a$ and $b$ being complex valued functions) generalizing the usual
analytic powers $\left\{  (z-z_{0})^{n}\right\}  _{n=0}^{\infty}$ in the sense
that locally when $z\rightarrow z_{0}$ they behave asymptotically like the
usual powers and under some additional conditions on the coefficients $a$ and
$b$ they form a complete system in the space of all solutions of the Vekua
equation in the same sense as the analytic powers $\left\{  (z-z_{0}%
)^{n}\right\}  _{n=0}^{\infty}$ form a complete system in the space of
analytic functions. Generalizations of the extension theorem, the Runge
theorem and of other important results about the convergence of corresponding
series are valid. The construction of formal powers is one of the main
problems of pseudoanalytic function theory. Recently it was solved
\cite{KrRecentDevelopments}, \cite{APFT} for a wide class of Vekua equations
of the form (\ref{Vekuamain}) which as was shown in the preceding section are
of main interest for studying problems for second order equations of the form
(\ref{maineq}).

The main ingredient for obtaining the explicit form of formal powers for a
certain Vekua equation is the generating sequence, a concept introduced by
Bers. If one knows a generating sequence for a given Vekua equation then the
construction of formal powers reduces to a simple algorithm. Here we briefly
explain the main ideas and steps refering the reader to \cite{Berskniga} and
\cite{APFT} for further details.

\subsection{Generating pair and generating sequence}

\begin{definition}
A pair of solutions $F$ and $G$ of a Vekua equation (\ref{VekuaGeneral}) in
$\Omega$ possessing partial derivatives with respect to the real variables $x$
and $y$ is said to be a generating pair if it satisfies the inequality
\begin{equation}
\operatorname{Im}(\overline{F}G)>0\qquad\text{in }\Omega. \label{GenPairCond}%
\end{equation}

\end{definition}

Condition (\ref{GenPairCond}) implies that every complex function $W$ defined
in a subdomain of $\Omega$ admits the unique representation $W=\phi F+\psi G$
where the functions $\phi$ and $\psi$ are real valued. Thus, the pair $(F,G)$
generalizes the pair $(1,i)$ which corresponds to usual complex analytic
function theory. The following expressions are known as characteristic
coefficients of the pair $(F,G)$%

\begin{align*}
a_{(F,G)}  &  =-\frac{\overline{F}G_{\overline{z}}-F_{\overline{z}}%
\overline{G}}{F\overline{G}-\overline{F}G},\qquad b_{(F,G)}=\frac
{FG_{\overline{z}}-F_{\overline{z}}G}{F\overline{G}-\overline{F}G},\\
A_{(F,G)}  &  =-\frac{\overline{F}G_{z}-F_{z}\overline{G}}{F\overline
{G}-\overline{F}G},\qquad B_{(F,G)}=\frac{FG_{z}-F_{z}G}{F\overline
{G}-\overline{F}G}.
\end{align*}

If $(F,G)$ is a generating pair of a Vekua equation (\ref{VekuaGeneral}) then
$a_{(F,G)}=a$ and $b_{(F,G)}=b$. \ The other two characteristic coefficients
are related to the concept of a derivative \cite{Berskniga}. The
$(F,G)$-derivative $\overset{\cdot}{W}=\frac{d_{(F,G)}W}{dz}$ of a
continuously differentiable function $W$ exists and has the form
\begin{equation}
\overset{\cdot}{W}=W_{z}-A_{(F,G)}W-B_{(F,G)}\overline{W} \label{FGder}%
\end{equation}
if and only if%
\[
W_{\overline{z}}=a_{(F,G)}W+b_{(F,G)}\overline{W}.
\]
Solutions of this equation are called $(F,G)$-pseudoanalytic functions.

\begin{definition}
\label{DefSuccessor}Let $(F,G)$ and $(F_{1},G_{1})$ - be two generating pairs
in $\Omega$. $(F_{1},G_{1})$ is called \ successor of $(F,G)$ and $(F,G)$ is
called predecessor of $(F_{1},G_{1})$ if%
\begin{equation}
a_{(F_{1},G_{1})}=a_{(F,G)}\qquad\text{and}\qquad b_{(F_{1},G_{1})}%
=-B_{(F,G)}\text{.} \label{charcoeffsuccessor}%
\end{equation}

\end{definition}

This definition arises naturally in relation to the notion of the
$(F,G)$-derivative due to the following fact.

\begin{theorem}
\label{ThBersDer}Let $W$ be an $(F,G)$-pseudoanalytic function and let
$(F_{1},G_{1})$ be a successor of $(F,G)$. Then $\overset{\cdot}{W}$ is an
$(F_{1},G_{1})$-pseudoanalytic function.
\end{theorem}

Thus, to the difference of analytic functions whose derivatives are again
analytic, the $(F,G)$-derivatives of pseudoanalytic functions are in general
solutions of another Vekua equation with the coefficients given by
(\ref{charcoeffsuccessor}). Obviously this process of construction of new
Vekua equations associated with the previous ones via relations
(\ref{charcoeffsuccessor}) can be continued and we arrive at the following definition.

\begin{definition}
\label{DefSeq}A sequence of generating pairs $\left\{  (F_{m},G_{m})\right\}
$, $m=0,\pm1,\pm2,\ldots$, is called a generating sequence if $(F_{m+1}%
,G_{m+1})$ is a successor of $(F_{m},G_{m})$. If $(F_{0},G_{0})=(F,G)$, we say
that $(F,G)$ is embedded in $\left\{  (F_{m},G_{m})\right\}  $.
\end{definition}

\begin{definition}
A generating sequence $\left\{  (F_{m},G_{m})\right\}  $ is said to have
period $\mu>0$ if $(F_{m+\mu},G_{m+\mu})$ is equivalent to $(F_{m},G_{m})$
that is their characteristic coefficients coincide.
\end{definition}

We will need the following notation introduced by Bers. The $(F,G)$-integral
is defined as follows
\[
\int_{\Gamma}Wd_{(F,G)}z=F(z_{1})\operatorname{Re}\int_{\Gamma}\frac
{2\overline{G}}{F\overline{G}-\overline{F}G}Wdz-G(z_{1})\operatorname{Re}%
\int_{\Gamma}\frac{2\overline{F}}{F\overline{G}-\overline{F}G}Wdz
\]
where $\Gamma$ is a rectifiable curve leading from $z_{0}$ to $z_{1}$.

Let $W$ be an $(F,G)$-pseudoanalytic function. Using a generating sequence in
which $(F,G)$ is embedded we can define the higher derivatives of $W$ by the
recursion formula%
\[
W^{[0]}=W;\qquad W^{[m+1]}=\frac{d_{(F_{m},G_{m})}W^{[m]}}{dz},\quad
m=0,1,\ldots\text{.}%
\]

A generating sequence defines an infinite sequence of Vekua equations. If for
a given (original) Vekua equation we know not only a corresponding generating
pair but the whole generating sequence, that is a couple of exact and
independent solutions for each of the Vekua equations from the infinite
sequence of equations corresponding to the original one, we are able to
construct an infinite system of solutions of the original Vekua equation as is
shown in the next definition.

\begin{definition}
\label{DefFormalPower}The formal power $Z_{m}^{(0)}(a,z_{0};z)$ with center at
$z_{0}\in\Omega$, coefficient $a$ and exponent $0$ is defined as the linear
combination of the generators $F_{m}$, $G_{m}$ with real constant coefficients
$\lambda$, $\mu$ chosen so that $\lambda F_{m}(z_{0})+\mu G_{m}(z_{0})=a$. The
formal powers with exponents $n=1,2,\ldots$ are defined by the recursion
formula%
\begin{equation}
Z_{m}^{(n)}(a,z_{0};z)=n\int_{z_{0}}^{z}Z_{m+1}^{(n-1)}(a,z_{0};\zeta
)d_{(F_{m},G_{m})}\zeta. \label{recformula}%
\end{equation}

\end{definition}

This definition implies the following properties.

\begin{enumerate}
\item $Z_{m}^{(n)}(a,z_{0};z)$ is an $(F_{m},G_{m})$-pseudoanalytic function
of $z$.

\item If $a^{\prime}$ and $a^{\prime\prime}$ are real constants, then
$Z_{m}^{(n)}(a^{\prime}+ia^{\prime\prime},z_{0};z)=a^{\prime}Z_{m}%
^{(n)}(1,z_{0};z)+a^{\prime\prime}Z_{m}^{(n)}(i,z_{0};z).$

\item The formal powers satisfy the differential relations%
\[
\frac{d_{(F_{m},G_{m})}Z_{m}^{(n)}(a,z_{0};z)}{dz}=nZ_{m+1}^{(n-1)}%
(a,z_{0};z).
\]

\item The asymptotic formulas
\begin{equation}
Z_{m}^{(n)}(a,z_{0};z)\sim a(z-z_{0})^{n},\quad z\rightarrow z_{0}
\label{asymptformulas}%
\end{equation}
hold.
\end{enumerate}

Assume now that
\begin{equation}
W(z)=\sum_{n=0}^{\infty}Z^{(n)}(a_{n},z_{0};z) \label{series}%
\end{equation}
where the absence of the subindex $m$ means that all the formal powers
correspond to the same generating pair $(F,G),$ and the series converges
uniformly in some neighborhood of $z_{0}$. It can be shown that the uniform
limit of pseudoanalytic functions is pseudoanalytic, and that a uniformly
convergent series of $(F,G)$-pseudoanalytic functions can be $(F,G)$%
-differentiated term by term. Hence the function $W$ in (\ref{series}) is
$(F,G)$-pseudoanalytic and its $r$th derivative admits the expansion
\[
W^{[r]}(z)=\sum_{n=r}^{\infty}n(n-1)\cdots(n-r+1)Z_{r}^{(n-r)}(a_{n}%
,z_{0};z).
\]
From this the Taylor formulas for the coefficients are obtained%
\begin{equation}
a_{n}=\frac{W^{[n]}(z_{0})}{n!}. \label{Taylorcoef}%
\end{equation}

\begin{definition}
Let $W(z)$ be a given $(F,G)$-pseudoanalytic function defined for small values
of $\left\vert z-z_{0}\right\vert $. The series%
\begin{equation}
\sum_{n=0}^{\infty}Z^{(n)}(a_{n},z_{0};z) \label{Taylorseries}%
\end{equation}
with the coefficients given by (\ref{Taylorcoef}) is called the Taylor series
of $W$ at $z_{0}$, formed with formal powers.
\end{definition}

The Taylor series always represents the function asymptotically:%
\begin{equation}
W(z)-\sum_{n=0}^{N}Z^{(n)}(a_{n},z_{0};z)=O\left(  \left\vert z-z_{0}%
\right\vert ^{N+1}\right)  ,\quad z\rightarrow z_{0}, \label{asympt}%
\end{equation}
for all $N$. This implies (since a pseudoanalytic function can not have a zero
of arbitrarily high order without vanishing identically) that the sequence of
derivatives $\left\{  W^{[n]}(z_{0})\right\}  $ determines the function $W$ uniquely.

If the series (\ref{Taylorseries}) converges uniformly in a neighborhood of
$z_{0}$, it converges to the function $W$.

\subsection{Convergence theorems\label{SubsectConvergence}}

S. Agmon and L. Bers \cite{AgmonBers} and L. Bers developed a theory of
expansions in pseudoanalytic formal powers which in its generality is
presented in \cite{Berskniga}, \cite{BersFormalPowers}. We do need here the
general results concerning a general Vekua equation (\ref{VekuaGeneral}).
Fortunately the situation with the main Vekua equation (\ref{Vekuamain}) in a
bounded simply connected domain under quite natural conditions on the function
$f$ is much easier than in the general case, and we have the following
expansion theorem and Runge theorem \cite{KrRecentDevelopments}, \cite{APFT}.

\begin{theorem}
\label{ThConvPer} Let $D$ be a disk of a finite radius $R$ and center $z_{0}$,
and $f\in C^{1}(\overline{D})$ be positive in $\overline{D}$. Then any
solution $W$ of (\ref{Vekuamain}) \ in $D$ admits a unique normally convergent
expansion\footnote{Following \cite{Berskniga}, \cite{Dettman} we shall say
that a sequence of functions $W_{n}$ converges normally in a domain $\Omega$
if it converges uniformly on every bounded closed subdomain of $\Omega$.} of
the form $W(z)=\sum_{n=0}^{\infty}Z^{(n)}(a_{n},z_{0};z)$.
\end{theorem}

\begin{theorem}
\label{ThRunge} any solution $W$ of (\ref{Vekuamain}) defined in a simply
connected domain can be expanded into a normally convergent series of formal
polynomials (linear combinations of formal powers with positive exponents).
\end{theorem}

\begin{remark}
This theorem admits a direct generalization onto the case of a multiply
connected domain (see \cite{BersFormalPowers}).
\end{remark}

We mention here another important result obtained by Menke in \cite{Menke}
which gives a useful estimate for the rate of convergence of the series from
the preceding theorem in the case when $W$ is a H\"{o}lder continuous function
up to the boundary of the domain of interest.

\begin{theorem}
\label{ThMenke} Let $W$ be a pseudoanalytic function in a domain $\Omega$
bounded by a Jordan curve and satisfy the H\"{o}lder condition on
$\partial\Omega$ with the exponent $\alpha$ ($0<\alpha\leq1$). Then for any
$\varepsilon>0$ and any natural $n$ there exists a pseudopolynomial of order
$n$ satisfying the inequality
\[
\left\vert W(z)-P_{n}(z)\right\vert \leq\frac{\operatorname*{Const}}%
{n^{\alpha-\varepsilon}}\qquad\text{for any }z\in\overline{\Omega}%
\]
where the constant does not depend on $n$, but only on $\varepsilon$.
\end{theorem}

The following statements are direct corollaries of the relations established
in section \ref{SectFactorization} between pseudoanalytic functions (solutions
of (\ref{Vekuamain})) and solutions of second-order elliptic equations, and of
the convergence theorems formulated above. Here we assume the existence of a
positive solution $u_{0}$ of (\ref{maineq}) in the domain $\Omega$ and the
function $f$ in (\ref{Vekuamain}) to be defined by $f=p^{1/2}u_{0}$ and belong
to $C^{1}(\overline{\Omega})$.

\begin{definition}
Let $u(z)$ be a given solution of the equation (\ref{maineq}) defined for
small values of $\left\vert z-z_{0}\right\vert $, and let $W(z)$ be a solution
of (\ref{Vekuamain}) constructed according to theorem \ref{CorConjugate}, such
that $\operatorname*{Re}W=p^{1/2}u$. The series
\begin{equation}
p^{-1/2}(z)\sum_{n=0}^{\infty}\operatorname*{Re}Z^{(n)}(a_{n},z_{0};z)
\label{TaylorForU}%
\end{equation}
with the coefficients given by (\ref{Taylorcoef}) is called the Taylor series
of $u$ at $z_{0}$, formed with formal powers.
\end{definition}

\begin{theorem}
\cite{KrJPhys06}, \cite{APFT} \label{ThExpMainEq}Let $u(z)$ be a solution of
(\ref{maineq}) defined for $\left\vert z-z_{0}\right\vert <R$. Then it admits
a unique expansion of the form
\[
u(z)=p^{-1/2}(z)\sum_{n=0}^{\infty}\operatorname*{Re}Z^{(n)}(a_{n},z_{0};z)
\]
which converges normally for $\left\vert z-z_{0}\right\vert <R$.
\end{theorem}

\begin{theorem}
\label{ThRungeSchr}An arbitrary solution of (\ref{maineq}) defined in a simply
connected domain where there exists a positive particular solution $u_{0}$
such that $f=p^{1/2}u_{0}\in C^{1}(\overline{\Omega})$ can be expanded into a
normally convergent series of formal polynomials multiplied by $p^{-1/2}$.
\end{theorem}

More precisely the last theorem has the following meaning. Due to Property 2
of formal powers we have that $Z^{(n)}(a,z_{0};z)$ for any Taylor coefficient
$a$ can be expressed through $Z^{(n)}(1,z_{0};z)$ and $Z^{(n)}(i,z_{0};z)$.
Then due to theorem \ref{ThRunge} any solution $W$ of (\ref{Vekuamain}) can be
expanded into a normally convergent series of linear combinations of
$Z^{(n)}(1,z_{0};z)$ and $Z^{(n)}(i,z_{0};z)$. Consequently, any solution of
(\ref{maineq}) can be expanded into a normally convergent series of linear
combinations of real parts of $Z^{(n)}(1,z_{0};z)$ and $Z^{(n)}(i,z_{0};z)$
multiplied by $p^{-1/2}$.

Obviously, for solutions of (\ref{maineq}) the results on the interpolation
and on the degree of approximation like, e.g., theorem \ref{ThMenke} are also valid.

Let us stress that theorem \ref{ThRungeSchr} gives us the following result.
The functions
\begin{equation}
\left\{  p^{-1/2}(z)\operatorname*{Re}Z^{(n)}(1,z_{0};z),\quad p^{-1/2}%
(z)\operatorname*{Re}Z^{(n)}(i,z_{0};z)\right\}  _{n=0}^{\infty}
\label{complsystem}%
\end{equation}
represent a complete system of solutions of (\ref{maineq}) in the sense that
any solution of (\ref{maineq}) can be represented by a normally convergent
series formed by functions (\ref{complsystem}) in any simply connected domain
$\Omega$ where a positive solution of (\ref{maineq}) exists, and the rate of
convergence of the series can be estimated with the aid of theorem
\ref{ThMenke}.

\subsection{Explicit construction of generating sequences and formal
powers\label{SectTheMainVekua}}

The results of section \ref{SectFactorization} show us that the theory of the
elliptic equation
\[
(\operatorname{div}p\operatorname{grad}+q)u=0
\]
is closely related to equation (\ref{Vekuamain}):
\begin{equation}
W_{\overline{z}}=\frac{f_{\overline{z}}}{f}\overline{W}. \label{Vekuamain1}%
\end{equation}
It is interesting that for this equation we always know a generating pair.
Namely, it is easy to see that the functions $F=f\quad$and\quad$G=\frac{i}{f}$
\ satisfy (\ref{Vekuamain1}) together with the condition (\ref{GenPairCond}).
Then the corresponding characteristic coefficients $A_{(F,G)}$ and $B_{(F,G)}$
have the form%
\[
A_{(F,G)}=0,\quad\text{\quad}B_{(F,G)}=\frac{f_{z}}{f},
\]
and the $(F,G)$-derivative according to (\ref{FGder}) is defined as follows%
\[
\overset{\cdot}{W}=W_{z}-\frac{f_{z}}{f}\overline{W}=\left(  \partial
_{z}-\frac{f_{z}}{f}C\right)  W.
\]

Due to Theorem \ref{ThBersDer} we obtain the following statement.

\begin{proposition}
\label{PrDer} Let $W$ be a solution of (\ref{Vekuamain1}). Then its
$(F,G)$-derivative, the function $w=\overset{\cdot}{W}$ is a solution of the
equation $\left(  \partial_{\overline{z}}+\frac{f_{z}}{f}C\right)  w=0$.
\end{proposition}

In spite of having given a generating pair for (\ref{Vekuamain1}) in general
it is not known how to construct a corresponding generating sequence necessary
for calculating the system of formal powers. Nevertheless a recent result from
\cite{KrRecentDevelopments}, \cite{APFT} which we formulate in the following
statement gives an answer to this question in a quite general situation.

\begin{theorem}
\label{ThGenSeq} Let $F=S(s)T(t)$ and $G=\frac{i}{S(s)T(t)}$ where $S$ and $T$
are arbitrary differentiable nonvanishing real valued functions, $\Phi=s+it$
is an analytic function of the variable $z=x+iy$ in $\Omega$ such that
$\Phi_{z}$ is bounded and has no zeros in $\Omega$. Then the generating pair
$(F,G)$ is embedded in the generating sequence $(F_{m},G_{m})$, $m=0,\pm
1,\pm2,\ldots$ in $\Omega$ defined as follows
\[
F_{m}=\left(  \Phi_{z}\right)  ^{m}F\quad\text{and\quad}G_{m}=\left(  \Phi
_{z}\right)  ^{m}G\quad\text{for even }m
\]
and%
\[
F_{m}=\frac{\left(  \Phi_{z}\right)  ^{m}}{S^{2}}F\quad\text{and\quad}%
G_{m}=\left(  \Phi_{z}\right)  ^{m}S^{2}G\quad\text{for odd }m.
\]

\end{theorem}

In order to appreciate the generality of this construction let us remind that
orthogonal coordinate systems in a plane are obtained (see \cite{Madelung})
from Cartesian coordinates $x$, $y$ by means of the relation%
\[
s+it=\Phi(x+iy)
\]
where $\Phi$ is an arbitrary analytic function. Quite often a transition to
more general coordinates is useful
\[
\xi=\xi(s),\quad\eta=\eta(t).
\]
$\xi$ and $\eta$ preserve the property of orthogonality. To illustrate the
point, besides the obvious example of Cartesian coordinates which are
generated by the analytic function $z$ we give some other examples taken from
\cite{Madelung}.

\begin{example}
\textbf{Polar coordinates }%
\[
s+it=\ln(x+iy),
\]%
\begin{equation}
s=\ln\sqrt{x^{2}+y^{2}},\quad t=\arctan\frac{y}{x}. \label{trueorthogcoord}%
\end{equation}
Usually the following new coordinates are introduced%
\[
r=e^{s}=\sqrt{x^{2}+y^{2}},\quad\varphi=t=\arctan\frac{y}{x}.
\]

\end{example}

\begin{example}
\textbf{Parabolic coordinates }%
\[
\frac{s+it}{\sqrt{2}}=\sqrt{x+iy},
\]%
\[
s=\sqrt{r+x},\quad t=\sqrt{r-x}.
\]
More frequently the parabolic coordinates are introduced as follows%
\[
\xi=s^{2},\quad\eta=t^{2}.
\]

\end{example}

\begin{example}
\textbf{Elliptic coordinates }%
\[
s+it=\arcsin\frac{x+iy}{\alpha},
\]%
\[
\sin s=\frac{s_{1}-s_{2}}{2\alpha},\quad\cosh t=\frac{s_{1}+s_{2}}{2\alpha}%
\]
where $s_{1}=\sqrt{(x+\alpha)^{2}+y^{2}}$, $s_{2}=\sqrt{(x-\alpha)^{2}+y^{2}}%
$. The substitution
\[
\xi=\sin s,\quad\eta=\cosh t
\]
is frequently used.
\end{example}

\begin{example}
\textbf{Bipolar coordinates}%
\[
s+it=\ln\frac{\alpha+x+iy}{\alpha-x-iy},
\]%
\[
\tanh s=\frac{2\alpha x}{\alpha^{2}+x^{2}+y^{2}},\quad\tan t=\frac{2\alpha
y}{\alpha^{2}-x^{2}-y^{2}}.
\]
The following substitution is frequently used%
\[
\xi=e^{-s},\quad\eta=\pi-t.
\]

\end{example}

The last theorem opens the way for explicit construction of formal powers
corresponding to the main Vekua equation (\ref{Vekuamain1}) in the case when
$f$ has the form
\begin{equation}
f=S(s)T(t) \label{fST}%
\end{equation}
and hence for explicit construction of complete systems of solutions for
corresponding second-order elliptic equations admitting a particular solution
of this form.

\section{Description of the method\label{SectDescription}}

We consider boundary value problems of Dirichlet, Neumann or mixed type for
the elliptic equation of the form (\ref{maineq}) in a bounded, simply
connected domain $\Omega\subset\mathbb{R}^{2}$. The main assumption required
for the applicability of the method of pseudoanalytic formal powers (MPFP) is
the existence in $\overline{\Omega}$ of a positive solution $u_{0}$ such that
the function $f=p^{1/2}u_{0}\in C^{1}(\overline{\Omega})$ be representable in
a separable form (\ref{fST}) in an orthogonal coordinate system. Let us stress
that very often such a particular solution $u_{0}$ is readily available. The
simplest example of such situation is when $q\equiv0$ and $p$ is of the form
(\ref{fST}). For example, the cases when $p(x,y)=X(x)Y(y)$ or $p=p(\sqrt
{x^{2}+y^{2}})$ frequently occur in practice \cite{Demidenko}.

When the equation of the form
\begin{equation}
(-\Delta+q(y))u(x,y)+\lambda^{2}u(x,y)=0 \label{Schr0}%
\end{equation}
is considered, it is sufficient to obtain a particular solution for the
ordinary differential equation
\begin{equation}
(-\frac{d^{2}}{dy^{2}}+q(y))h(y)=0. \label{ordinary}%
\end{equation}
Then a particular solution of (\ref{Schr0}) can be constructed as follows
\begin{equation}
u_{0}(x,y)=e^{\lambda x}h(y). \label{u0}%
\end{equation}
It has a convenient separable form. Notice that in this example we come to an
important open problem. It is related to the requirement that $u_{0}$ should
be different from zero in the domain of interest. Meanwhile in many
practically significant situations it is easy to guarantee that $h(y)\neq0$
when $(x,y)\in\overline{\Omega}$, sometimes this condition becomes a
considerable obstacle. Moreover, when $\lambda$ in (\ref{Schr0}) is purely
imaginary, the solution (\ref{u0}) is not acceptable because it is no longer
real valued. In this case one should\ take instead of $e^{ikx}$, where
$\lambda=ik$, a solution in the form of $\sin kx$ or $\cos kx$ but then for a
big domain or large $k$ one cannot avoid the appearance of zeros of the
resulting particular solution $u_{0}$ and in this case the proposed scheme in
general does not work.

One possibility to overcome this problem is to include under consideration a
complex valued particular solution $u_{0}$ but then we would need to consider
a corresponding bicomplex main Vekua equation (see \cite{KrJPhys06} and
\cite{APFT}). Then the whole algorithm for the construction of generating
sequences and formal powers would go through with no modification compared to
the complex case, but up to now there is no proof of the completeness of the
system of formal powers for a bicomplex Vekua equation. As a consequence there
is no guarantee that the infinite system of exact solutions obtained similar
to (\ref{complsystem}) will be complete in the space of solutions of
(\ref{maineq}) in $\Omega$. Our conjecture is that at least in the case when
$u_{0}(x,y)=g(x)h(y)$ where $g$ and $h$ are complex valued nonvanishing
functions the system of formal powers for the corresponding bicomplex main
Vekua equation is complete in the same sense as was established earlier for
the complex case. We continue this discussion in section \ref{SectEigen} where
we use complex valued particular solutions of the form (\ref{u0}) for solving
eigenvalue problems for operators of the form $-\Delta+q(y)$.

Turning back to equation (\ref{maineq}) we assume that it admits a positive
solution $u_{0}$ in the domain $\Omega$ such that $f=p^{1/2}u_{0}\in
C^{1}(\overline{\Omega})$ is representable in a separable form (\ref{fST}) and
that $\Phi=s+it$ is an analytic function of the variable $z=x+iy$ in $\Omega$
such that $\Phi_{z}$ is bounded and has no zeros in $\Omega$. Then applying
theorem \ref{ThGenSeq} one can construct a corresponding generating sequence.
Construction of formal powers $\left\{  Z^{(n)}(1,z_{0};z),\quad
Z^{(n)}(i,z_{0};z)\right\}  _{n=0}^{\infty}$ reduces then to the recursive
algorithm described in Definition \ref{DefFormalPower}, and in this way one
obtains the complete system of solutions for (\ref{maineq}) in $\Omega$ given
by (\ref{complsystem}). By construction $\operatorname*{Re}Z^{(0)}%
(i,z_{0};z)\equiv0$. Taking this into account we introduce the notations
\begin{align*}
u_{1}(z) &  =p^{-1/2}(z)\operatorname*{Re}Z^{(1)}(1,z_{0};z),\quad
u_{2}(z)=p^{-1/2}(z)\operatorname*{Re}Z^{(1)}(i,z_{0};z),\\
u_{3}(z) &  =p^{-1/2}(z)\operatorname*{Re}Z^{(2)}(1,z_{0};z),\quad
u_{4}(z)=p^{-1/2}(z)\operatorname*{Re}Z^{(2)}(i,z_{0};z),\ldots
\end{align*}
and obtain the complete system of solutions for (\ref{maineq}) given by
$\left\{  u_{0},u_{1},u_{2},\ldots\right\}  $. We look for an approximate
solution of a boundary value problem for (\ref{maineq}) in the form
\begin{equation}
u^{N}=\sum_{k=0}^{N}b_{k}u_{k}\label{uN}%
\end{equation}
where $b_{k}$ are real coefficients which should be found from boundary
conditions. To obtain $N+1$ equations for finding $\left\{  b_{k}\right\}
_{k=0}^{N}$ one can use, e.g., the collocation method. Chosing $N+1$ points
$\zeta_{j}\in\partial\Omega$ we obtain $N+1$ equations
\[
\sum_{k=0}^{N}b_{k}B[u_{k}](\zeta_{j})=v(\zeta_{j}),\quad j=\overline{0,N}%
\]
where $v$ is a given function and $B$ is the linear operator of the boundary
condition. For the Dirichlet condition one has $B[u]=u$ and for the Neumann
condition, $B[u]=\frac{\partial u}{\partial\overrightarrow{n}}$ -- the normal
derivative of $u$. Finding $\left\{  b_{k}\right\}  _{k=0}^{N}$ we have the
approximate solution $u^{N}$.

Thus, the proposed here MPFP belongs to the class of boundary methods because
due to the linearity of the problem the function (\ref{uN}) is an exact
solution of (\ref{maineq}) in $\Omega$ and only boundary conditions should be
approximated. An estimate for the rate of convergence of the method is given
in theorem \ref{ThMenke}. In the next section we discuss the numerical
realization of MPFP and results of numerical tests.

\section{Approximate solution of boundary value problems}

As a first example we considered the Dirichlet problem for the equation
\begin{equation}
\left(  -\Delta+c^{2}\right)  u=0\label{Yukawa}%
\end{equation}
where $c$ is a real constant. The interest in this relatively simple equation
is not due to its numerical simplicity. In fact this is not the case,-
numerical solution of this equation is not less difficult than that of an
equation with $c$ being a reasonably good function with a range of values
comparable with $c$. The attractiveness of this example consists in the
possibility to calculate\ a large number of the functions $u_{k}$ (see the
preceding section) symbolically, using an appropriate software for symbolic
calculations like Mathematica (Wolfram), Maple or Matlab. In this work we used
Matlab 2006 and a PC of 2 GB in RAM and a processor of 1.73 GHz.
Implementation of the symbolically calculated base functions $u_{k}$ gives us
the possibility to estimate the accuracy of the MPFP itself without
considering the precision of recursive numerical integrations. We also compare
the results obtained using symbolically calculated $u_{k}$ with the results
obtained purely numerically.

For equation (\ref{Yukawa}) it is easy to propose a positive particular
solution. It can be chosen, e.g., as $f=e^{cy}$. Then the first functions
$u_{k}$ constructed as described in the preceding section taking as a center
of the formal powers the origin will have the form \cite{APFT}%
\begin{align}
u_{0}(x,y) &  =e^{cy},\qquad u_{1}(x,y)=xe^{cy},\qquad u_{2}(x,y)=-\frac
{\sinh(cy)}{c},\nonumber\\
& \label{uk}\\
u_{3}(x,y) &  =\left(  x^{2}-\frac{y}{c}\right)  e^{cy}+\frac{\sinh(cy)}%
{c^{2}},\qquad u_{4}(x,y)=-\frac{2x\sinh(cy)}{c},\ldots\text{.}\nonumber
\end{align}
It is interesting to mention that using Matlab we obtained the first 101
functions of this system calculated symbolically. According to theorem
\ref{ThRungeSchr} this system of solutions is complete in any bounded simply
connected domain containing the origin. First we show results obtained with
the help of the system of functions $u_{k}$ calculated symbolically.

\subsection{Numerical results obtained with symbolically calculated base
functions}

We begin with the unitary disk $D$ with center in the origin. As a test exact
solution we take the function
\begin{equation}
u=e^{cx}. \label{exact}%
\end{equation}
Thus, the problem we consider is to solve (\ref{Yukawa}) in $D$ with the
boundary condition $\left.  u\right\vert _{\partial D}=e^{cx}$. We look for an
approximate solution $u^{N}$ in the form (\ref{uN}) with the base functions
(\ref{uk}). We use the collocation method for satisfying the boundary
condition, the collocation points are distributed uniformly on $\partial D$.
Their number is equal to the number of solutions $u_{k}$.

According to theory from subsection \ref{SubsectConvergence} the coefficients
$b_{k}$ in (\ref{uN}) are obtained in the case under consideration from the
Taylor coefficients which appear in (\ref{TaylorForU}). More precisely we have
that according to theorem \ref{ThExpMainEq} the solution $u$ can be
represented as follows
\[
u(z)=\sum_{n=0}^{\infty}\operatorname*{Re}Z^{(n)}(a_{n},0;z)=\sum
_{n=0}^{\infty}\left(  a_{n}^{\prime}\operatorname*{Re}Z^{(n)}(1,0;z)+a_{n}%
^{\prime\prime}\operatorname*{Re}Z^{(n)}(i,0;z)\right)
\]
where $a_{n}=a_{n}^{\prime}+ia_{n}^{\prime\prime}$ are the Taylor coefficients
given by (\ref{Taylorcoef}). In the case of the exact solution (\ref{exact})
the Taylor coefficients have the form \cite[Sect. 7.3]{APFT}%
\[
a_{n}=\frac{c^{n}}{n!}(1+i).
\]
Thus, the exact values for the coefficients $b_{k}$ from (\ref{uN}) in our
example are as follows%
\[
b_{0}=1,\quad b_{1}=b_{2}=c,\quad b_{3}=b_{4}=\frac{c^{2}}{2},\ldots.
\]

Having compared the numerically calculated constants $b_{k}$ which we denote
by $\widetilde{b}_{k}$ for $N=34$ with their exact values in the case $c=1$ we
obtained their coincidence up to $10^{-14}$ for every $k=0,\ldots,34$. For
smaller values of $c$ the situation is the same. The difference between
$\widetilde{b}_{k}$ and $b_{k}$ tends to become larger for larger values of
$c$. In Table \ref{Table1} we show results for $c=5$ and $N=34$.

\begin{center}
Table 1.\label{Table1} Comparison of the values of $\widetilde{b}_{k}$ and
$b_{k}$ as $k$ increases%

\begin{tabular}
[c]{|c|c|c|c|}\hline
$k$ & The values of $\widetilde{b}_{k}$ & The values of $b_{k}$ & $\left\vert
\widetilde{b}_{k}-b_{k}\right\vert $\\\hline
$5$ & $20.83333333333382$ & $20.83333333333333$ & $0.00000000000048$\\\hline
$8$ & $26.04166666666448$ & $26.04166666666667$ & $0.00000000000219$\\\hline
$13$ & $15.50099206509864$ & $15.50099206349206$ & $0.00000000160657$\\\hline
$17$ & $5.38228885848900$ & $5.38228891093474$ & $0.00000005244574$\\\hline
$25$ & $0.19601580023149$ & $0.19603324996120$ & $0.00001744972971$\\\hline
$31$ & $0.00743227875004$ & $0.00729290364439$ & $0.00013937510565$\\\hline
$34$ & $0.00172611010091$ & $0.00214497166011$ & $0.00041886155920$\\\hline
\end{tabular}

\end{center}

The maximum number of functions $u_{k}$ that we used here is limited not by
the possibility of obtaining them symbolically but rather by the time required
for numerical calculations involving the corresponding quite long symbolic expressions.

In the following two tables, the convergence of MPFP is shown by comparison of
the maximum absolute error obtained for different values of $N$, for the case
$c=1$ and $c=5$.

\begin{center}
Tables 2 and 3\label{Table2} Maximum absolute error depending on $N$ for $c=1$
and $c=5$

\begin{tabular}
[c]{|c|c|}\hline
$N$ & Maximum absolute error\\\hline
$8$ & $0.00698626935341$\\\hline
$14$ & $2.534633673767495\times10^{-5}$\\\hline
$22$ & $1.432881036045330\times10^{-9}$\\\hline
$28$ & $4.276579090856103\times10^{-13}$\\\hline
$32$ & $1.776356839400251\times10^{-15}$\\\hline
$36$ & $8.881784197001252\times10^{-16}$\\\hline
$38$ & $1.110223024625157\times10^{-15}$\\\hline
\end{tabular}
\qquad%
\begin{tabular}
[c]{|c|c|}\hline
$N$ & Maximum absolute error\\\hline
$6$ & $3.59578971016677\times10^{2}$\\\hline
$14$ & $22.38029523897584$\\\hline
$22$ & $0.73431266884919$\\\hline
$32$ & $0.00194275813006$\\\hline
$44$ & $0.59167057031573\times10^{-7}$\\\hline
$54$ & $0.72795103278622\times10^{-11}$\\\hline
$60$ & $8.781864124784988\times10^{-14}$\\\hline
\end{tabular}

\end{center}

In Table 4 for a fixed number ($N=34$) of the base functions $u_{k}$ we show
the dependence of the maximum absolute error on the parameter $c$. Here we
also indicate the maximum absolute error obtained for the same problem using
the standard PDE tool of Matlab.

\begin{center}
Table 4.\label{Table4} Performance of MPFP compared to Matlab's PDE tool in
terms of the maximum absolute error for increasing values of $c$ as $N=34$%

\begin{tabular}
[c]{|c|c|c|}\hline
$c$ & Maximum absolute error of MPFP & PDE tool (2129 nodes)\\\hline
$0.1$ & $0.89\times10^{-15}$ & $1.5\times10^{-6}$\\\hline
$0.5$ & $0.26\times10^{-14}$ & $4.5\times10^{-6}$\\\hline
$1$ & $0.12\times10^{-14}$ & $1.6\times10^{-4}$\\\hline
$2$ & $0.14\times10^{-10}$ & $1.4\times10^{-3}$\\\hline
$5$ & $0.29\times10^{-3}$ & $3.0\times10^{-2}$\\\hline
$10$ & $4.06\times10^{2}$ & $8.0$\\\hline
\end{tabular}

\end{center}

As it can be observed in the last table the result of application of MPFP in
the case of $c=10$ is less satisfactory as that of PDE tool. This is due to
the fact that for larger values of $c$ one should consider a bigger $N$. In
Table 5 we show the absolute error of MPFP for $c=10$ and $N\geq42$.

\begin{center}
Table 5. Improvement in the maximum absolute error due to MPFP as the number
of functions $u_{k}$ keeps increasing%

\begin{tabular}
[c]{|c|c|}\hline
$N$ & Maximum absolute error of MPFP\\\hline
$42$ & $3.89$\\\hline
$44$ & $3.25$\\\hline
$46$ & $1.81$\\\hline
$48$ & $0.81$\\\hline
$50$ & $0.41$\\\hline
$52$ & $0.10$\\\hline
\end{tabular}

\end{center}

Thus, one can see that for $N\geq42$ the result obtained with the aid of MPFP
is more accurate than that given by Matlab. We stress that in the case of
using MPFP a system of $N+1$ linear algebraic equations is solved which means
solution of dozens of equations instead of thousands required by the finite
element method implemented in the PDE tool.

We experimented also with the shape of the domain. We considered the elliptic
form as well as a unitary disk with a triangle shaped deformation. In the
first case it is possible to see how the maximum absolute error increases with
the excentricity $e$, Table 6. Here in all cases the area of the considered
ellipses was kept constant, equal to $\pi$, while the excentricity was being increased.

\begin{center}
Table 6. Maximum absolute error for different values of the excentricity of
the elliptic domain with the area of the domain being equal to $\pi$. The case
$e=0$ corresponds to the unitary disc.%

\begin{tabular}
[c]{|c|c|c|c|c|c|c|}\hline
$N$ & $e=0$ & $e=0.5$ & $e=0.7$ & $e=0.9$ & $e=0.95$ & $e=0.99$\\\hline
$30$ & $2.2\times10^{-14}$ & $0.4\times10^{-13}$ & $0.5\times10^{-13}$ &
$0.3\times10^{-12}$ & $0.2\times10^{-11}$ & $1\times10^{-10}$\\\hline
\end{tabular}

\end{center}

For the case of the domain with a triangular deformation (see Fig. 1),
%TCIMACRO{\FRAME{ftbpFU}{4.0171in}{3.0173in}{0pt}{\Qcb{The unitary disk with a
%triangular deformation.}}{}{fig1.jpg}{\special{ language "Scientific Word";
%type "GRAPHIC";  maintain-aspect-ratio TRUE;  display "USEDEF";
%valid_file "F";  width 4.0171in;  height 3.0173in;  depth 0pt;
%original-width 5.8332in;  original-height 4.3751in;  cropleft "0";
%filename 'Fig1.jpg';file-properties "XNPEU";}}}%
%BeginExpansion
\begin{figure}
[ptb]
\begin{center}
\includegraphics[
natheight=4.375100in,
natwidth=5.833200in,
height=3.0173in,
width=4.0171in
]%
{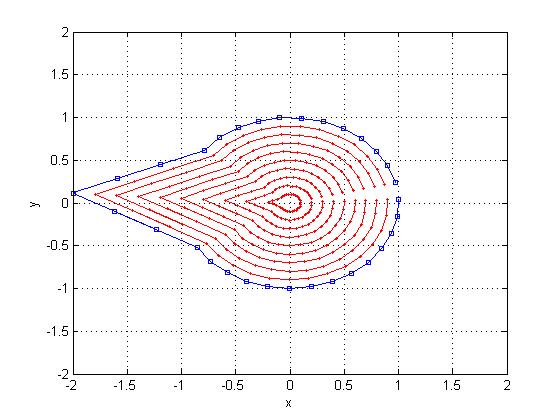}%
\caption{The unitary disk with a triangular deformation.}%
\end{center}
\end{figure}
%EndExpansion
the errors were tested for different heights of the peak and different values
of $c$ and $N$ with satisfactory results. In Table 7 we present the maximum
absolute error of the approximate solution of the boundary value problem in
dependence on the height of the triangular peak over the unitary circunference.

\begin{center}
Table 7. Maximum absolute error for $N=31$, $c=1$ and different heights of the peak%

\begin{tabular}
[c]{|c|c|}\hline
Height of the peak over the unitary disk & Maximum absolute error\\\hline
$0.5$ & $0.92\times10^{-12}$\\\hline
$0.7$ & $0.62\times10^{-11}$\\\hline
$1.0$ & $0.72\times10^{-10}$\\\hline
\end{tabular}
\bigskip
\end{center}

\subsection{Results obtained with numerically calculated base functions}

The use of the numerically calculated base functions which we denote by
$\widetilde{u}_{k}$ poses the natural question about the accuracy of their
calculation. Consideration of equation (\ref{Yukawa}) gives us the possibility
to compare $\widetilde{u}_{k}$ with the symbolically calculated exact
solutions (\ref{uk}). In the following table we give the difference between
$u_{k}$ and $\widetilde{u}_{k}$ for  $c=1$.

\begin{center}
Table 8. Maximum absolute error of calculation of the base functions for $c=1$%

\begin{tabular}
[c]{|c|c|}\hline
$k$ & $\left\vert u_{k}-\widetilde{u}_{k}\right\vert $\\\hline
$1$ & $0.00000667646050\times10^{-4}$\\\hline
$5$ & $0.00022726096338\times10^{-4}$\\\hline
$11$ & $0.01492959925020\times10^{-4}$\\\hline
$16$ & $0.07644765777970\times10^{-4}$\\\hline
$20$ & $0.22545767650040\times10^{-4}$\\\hline
\end{tabular}

\end{center}

Regarding the results given in the last table it is important to note that in
fact the weight of $u_{k}$ in the expansion of a solution decreases as
$\left(  \frac{k+1}{2}\right)  !$ or $\left(  \frac{k}{2}\right)  !$ for odd
or even $k$ respectively. This is due to the factor $1/n!$ in the definition
of the Taylor coefficients (\ref{Taylorcoef}). That is in fact the real
accuracy of calculation of the base functions would be given by $\left\vert
u_{k}-\widetilde{u}_{k}\right\vert /K!$ where $K=\left\{
\begin{array}
[c]{c}%
\left(  \frac{k+1}{2}\right)  !\text{ for }k\text{ odd}\\
\left(  \frac{k}{2}\right)  !\text{ for }k\text{ even}%
\end{array}
\right.  $. It is easy to see that in the case of the results given in Table 7
one then obtains $\left\vert u_{20}-\widetilde{u}_{20}\right\vert
/10!\simeq6.213\times10^{-12}$ which is a remarkably good agreement.

For the numerical computation of $\widetilde{u}_{k}$ we implemented the
following procedure. Before integrating on each new step according to
(\ref{recformula}) along segments joining the center of the formal powers with
points on the boundary of the domain the integrand was represented as a cubic
spline which then was integrated using the standard Matlab routine for
integration of splines. This procedure is simple but clearly not optimal.
Nevertheless the approximate results presented in this work show that even
such integration procedure gives satisfactory agreement between the exact base
functions and those calculated numerically.

The accuracy of the approximate solution obtained with the aid of
$\widetilde{u}_{k}$ in our numerical tests did not differ significantly from
that of the solution obtained using the exactly calculated $u_{k}$. The order
of the maximum absolute error for a given $N$ and $c$ coincided in both cases.
Hence here we present results corresponding to another test problem for which
we did not have the exactly calculated base functions.

Consider the equation
\begin{equation}
\left(  -\Delta+\frac{e^{y}}{4}\right)  u(x,y)=0.\label{secondexample}%
\end{equation}
An exact solution for this equation can be found using the fact that the
change of variables $\xi=e^{\frac{y}{2}}\cos\frac{x}{2}$, $\eta=e^{\frac{y}%
{2}}\sin\frac{x}{2}$ leads to (\ref{Yukawa}) in the new variables. Thus, e.g.,
the function $u(x,y)=\exp(e^{\frac{y}{2}}\cos\frac{x}{2})$ is an exact
solution of (\ref{secondexample}). Consequently, as a test problem we can
consider the problem of finding a solution of (\ref{secondexample}) in
$\Omega$, satisfying the boundary condition
\begin{equation}
u(x,y)=\exp(e^{\frac{y}{2}}\cos\frac{x}{2}),\quad(x,y)\in\partial
\Omega.\label{test2}%
\end{equation}
In order to construct a particular solution $u_{0}$ in a separable form we
solve numerically the ordinary differential equation
\[
\left(  -\frac{d^{2}}{dy^{2}}+\frac{e^{y}}{4}\right)  u_{0}(y)=0.
\]
The obtained solution we then use for constructing the system of functions
$\left\{  u_{k}\right\}  $. Some results on the accuracy of the approximate
solution are given in the following table.

\begin{center}
\bigskip Table 9. Maximum absolute error of the approximate solution of the
test problem (\ref{secondexample}), (\ref{test2}) considered in a unitary disk
in dependence on $N$%

\begin{tabular}
[c]{|c|c|}\hline
$N$ & Maximum absolute error\\\hline
$4$ & $0.079$\\\hline
$6$ & $0.021$\\\hline
$8$ & $0.005$\\\hline
$10$ & $0.001$\\\hline
$12$ & $0.0004$\\\hline
$14$ & $0.000099$\\\hline
$16$ & $0.000020$\\\hline
$18$ & $0.0000033$\\\hline
$20$ & $0.00000060$\\\hline
$28$ & $0.00000000072$\\\hline
$32$ & $0.00000000028$\\\hline
\end{tabular}

\end{center}

\section{Approximate solution of eigenvalue problems\label{SectEigen}}

In this section we consider the application of MPFP to solution of eigenvalue
problems for operators of the form $-\Delta+q(y)$. \ For simplicity we keep
working with the Dirichlet boundary conditions and suppose that $q$ is
continuous and $q(y)\geq0$ in $\overline{\Omega}$. Then the spectrum of the
operator is discrete and positive. As was explained in section
\ref{SectDescription} for the equation%
\begin{equation}
(-\Delta+q(y))u_{0}(x,y)=\lambda^{2}u_{0}(x,y)\label{Schrodwitheigen}%
\end{equation}
it is easy to propose a particular solution in a separable form for any value
of $\lambda$. We are interested here in positive values, and hence a natural
choice of a nonvanishing solution would be $u_{0}(x,y)=e^{i\lambda x}h(y)$
where $h(y)$ is a positive solution of (\ref{ordinary}). As was observed in
section \ref{SectDescription} the completeness of the system of solutions
$\left\{  u_{k}\right\}  _{k=0}^{\infty}$ obtained in this case is up to now
an open problem due to the fact that $u_{0}$ is complex valued and one should
consider bicomplex pseudoanalytic formal powers for which the whole theory is
still underdeveloped. Nevertheless we used the constructed system of exact
solutions $\left\{  u_{k}\right\}  _{k=0}^{\infty}$ for finding the
eigenvalues $\lambda^{2}$ in the following way. Assuming that $\left\{
u_{k}\right\}  _{k=0}^{\infty}$ is complete in the same sense as was proved in
the case of the real-valued particular solution $u_{0}$ (subsection
\ref{SubsectConvergence}) we have then that if a nontrivial solution $u$ of
(\ref{Schrodwitheigen}) exists satisfying the boundary condition $\left.
u\right\vert _{\partial\Omega}=0$ then $u\simeq\sum_{k=0}^{N}b_{k}u_{k}$ and
the coefficients $b_{k}$ are such that the trivial boundary condition is
approximately fulfilled. This means that one can require that $\sum_{k=0}%
^{N}b_{k}u_{k}(z_{j})=0$ for $z_{j}\in\partial\Omega$ and $j=\overline{0,N}$.
This is possible iff the determinant of the matrix $U=\left(  u_{jk}\right)
_{j,k=0}^{N}$ vanishes where $u_{jk}$ $=u_{k}(z_{j})$. The determinant of $U$
for a fixed $N$ is a function of $\lambda$. Thus, the problem of finding
eigenvalues reduces to the problem of finding zeros of the function $\det
U(\lambda)$.

As a test problem we considered the problem of calculating the eigenvalues of
the Dirichlet problem for the Helmholtz equation $\left(  \Delta+\lambda
^{2}\right)  u=0$. For every $\lambda$ \ a system of exact solutions $\left\{
u_{k}\right\}  _{k=0}^{\infty}$ can be constructed using (\ref{uk}) where $c$
should be replaced by $i\lambda$. Then following the described scheme we
looked for zeros of $\det U(\lambda)$. As it is well known (see, e.g.,
\cite{Coleman}) the eigenvalues of the Dirichlet problem for the Helmholtz
equation in a unitary disk are squares of zeros of Bessel functions $J_{n}%
(x)$. Our numerical experiments showed that a relatively small value of $N=21$
was needed for computing the first five eigenvalues with the accuracy of four
decimals. With $N=23$ we obtained six first eigenvalues with the same
accuracy. Thus, indeed, the MPFP is clearly competitive in solving eigenvalue
problems for elliptic operators. In applying the MPFP we detected a similar
problem to that described by Alexidze in \cite[Sect. 1.13]{Alexidze} where the
method of fundamental solutions (or auxiliary sources) was applied to
eigenvalue problems. The considered determinant shows a very fast decrement
(in spite of this the method gives good numerical results). \cite{Alexidze}
contains references to other publications where different ways of using the
knowledge of a system of exact solutions for numerical solution of eigenvalue
problems were studied. In this direction further research is needed.

\section{Conclusions}

A new approach for solving boundary value and eigenvalue problems for elliptic
operators in bounded planar domains is proposed. It is based on some classical
and some new results from pseudoanalytic function theory which allow one to
construct complete systems of solutions of the elliptic equations. We showed
the practical applicability of the numerical method based on this
construction, studied the rate of its convergence, accuracy and other
parameters of its performance.

\end{document}